\input{STU.STY}
\usepackage{xcolor}
\newfont{\cyr}{wncyr8}
\newfont{\cyb}{wncyr8}
\newtheorem{thm}{Theorem}[section]
\newcommand{\bthm}{\begin{thm}}
\newcommand{\ethm}{\end{thm}}

\newtheorem{prop}[thm]{Proposition}
\newcommand{\bprp}{\begin{prop}}
\newcommand{\eprp}{\end{prop}}

\newtheorem{fact}[thm]{Fact}
\newcommand{\bfct}{\begin{fact}}
\newcommand{\efct}{\end{fact}}

\newtheorem{prob}[thm]{Problem}
\newcommand{\bprb}{\begin{prob}}
\newcommand{\eprb}{\end{prob}}

\newtheorem{quest}[thm]{Question}
\newcommand{\bqtn}{\begin{quest}}
\newcommand{\eqtn}{\end{quest}}

\newtheorem{lem}[thm]{Lemma}
\newcommand{\blem}{\begin{lem}}
\newcommand{\elem}{\end{lem}}

\newtheorem{claim}[thm]{Claim}
\newcommand{\bclm}{\begin{claim}}
\newcommand{\eclm}{\end{claim}}

\newtheorem{cor}[thm]{Corollary}
\newcommand{\bcor}{\begin{cor}}
\newcommand{\ecor}{\end{cor}}

\newtheorem{conj}[thm]{Conjecture}
\newcommand{\bcnj}{\begin{conj}}
\newcommand{\ecnj}{\end{conj}}

\theoremstyle{definition}
\newtheorem{defn}[thm]{Definition}
\newcommand{\bdfn}{\begin{defn}}
\newcommand{\edfn}{\end{defn}}

\newtheorem{spec}[thm]{Specializing}
\newcommand{\bspc}{\begin{spec}}
\newcommand{\espc}{\end{spec}}

\theoremstyle{remark}
\newtheorem{rem}[thm]{Remark}
\newcommand{\brem}{\begin{rem}}
\newcommand{\erem}{\end{rem}}

\newtheorem{cnv}[thm]{Convention}
\newcommand{\bcnv}{\begin{cnv}}
\newcommand{\ecnv}{\end{cnv}}

\newtheorem{exam}[thm]{Example}
\newcommand{\bexm}{\begin{exam}}
\newcommand{\eexm}{\end{exam}}

\newtheorem{exercise}[thm]{Exercise}
\newcommand{\bexr}{\begin{exercise}}
\newcommand{\eexr}{\end{exercise}}

\newtheorem{thmy}{\textbf{Theorem}}

\usepackage{pdfpages}
\usepackage{amsxtra}

\begin{document}
\label{begin-art}

\title{The Baire property and precompact duality}

\date{\today}

\author[M. Ferrer, S. Hern\'andez, I. Sep\'ulveda, F. J. Trigos-Arrieta]
{M. Ferrer, S. Hern\'andez, I. Sep\'ulveda, F. J. Trigos-Arrieta}
\address{Universitat Jaume I, Departamento de Matem\'{a}ticas,
Campus de Riu Sec, 12071 Castell\'{o}n, Spain.}
\email{mferrer@mat.uji.es, hernande@mat.uji.es, isepulve@mat.uji.es}
\address{Departament of Mathematics,
California State University, Bakersfield, Bakersfield California, USA.}
\email{jtrigos@csub.edu}

\thanks{The first three-listed authors acknowledge partial support by the Spanish Ministerio de Econom\'{i}a y Competitividad, grant: MTM/PID2019-106529GB-I00 (AEI/FEDER, EU) and by the Universitat Jaume I, grant UJI-B2022-39. The fourth author is very grateful to the Departament de Mathem\`atiques of Universitat Jaume I {and the Institut de Matem\`atiques i Aplicacions de Castell\'o (IMAC)} for generous support during his sabbatical stay in Fall 2023.}

\begin{abstract}
We prove that if $G$ is a totally bounded abelian group \st\ its dual group $\widehat{G}_p$ equipped with
the finite-open topology is a Baire group, then every compact subset of $G$ must be finite.
This solves an open question by Chasco, Dom\'inguez and Tkachenko. {Among other consequences, we obtain
an example of a group that is $g$-dense in its completion but is not $g$-barrelled. This solves a question
proposed by Au\ss enhofer and Dikranjan.}
\end{abstract}



\maketitle \setlength{\baselineskip}{24pt}

\section{Introduction and basic facts}
Let $G$ be a \tb\ abelian group and let $\widehat{G}_p$ denote its dual group equipped with the finite-open topology $\sigma(\widehat{G},G)$
(of poinwise convergence on the elements of $G$). We investigate here the consequences on the group $G$ of imposing the Baire property on the dual group $\widehat{G}_p$.
In particular, we focus on the question, set by Chasco, Dom\'inguez and Tkachenko in \cite{ChascoDominguezTkachenko} and repeated in \cite{FerHerTka,AusDik:2021}, of whether a totally bounded group $G$,
whose dual group $\widehat{G}_p$ is Baire, may contain infinite compact subsets. In fact, {the authors of \cite{ChascoDominguezTkachenko} prove that if $G$ is a bounded torsion group
and $\widehat{G}_p$ is Baire, then $G$ does not contain infinite compact subsets (subsequently, it was proved in \cite{AusDik:2021} that the constraint
``bounded torsion group" can be relaxed to  ``torsion group".}) {Here, we complete 
these re\-sults by proving that the same is true in general}. That is, if $G$ is a totally bounded group, then the Baire property of $\widehat{G}_p$ implies
that $G$ does not contain infinite compact subsets. As a consequence, several applications of this result are provided.
For example, {we have that any Baire totally bounded group $G$, {without infinite compact subsets},  must be automatically Pontryagin reflexive since its dual $\widehat{G}_p$ contains no infinite compact subsets,} a phenomenon which
contains most known examples of totally bounded reflexive groups. {Finally, we give an example of a group that is $g$-dense
in its completion but is not $g$-barrelled. This solves a question proposed by Au\ss en\-hofer and Dikranjan \cite{AusDik:2021}.}

Unless we say otherwise, all our groups are commutative {and our spaces Tychonoff.} 
With $\ZZ$ we denote the group of integers. {At times, we consider $\ZZ$ as a topological group, while at other times, we regard it as the underlying group without any topology.}
Si\-mi\-lar\-ly, by $\ZZ_q$ we denote the cyclic group of order $q \in \NN$. Set $\omega:=\{0,1,2,...\}$. If $S$ is a set, we denote by $S^{<\omega}$ the set of finite subsets of $S$.
Our model for $\TT$ is the group $([0,1),+$~mod~$1$) with the topology inherited from $\RR$ {and the points $0$ and $1$ identified}, unless we explicitly say otherwise
(We use the model $((-1/2,1/2],+$~mod~$1$) in the proof of {Theorem} \ref{teor1}),
when we need to see it as a \tg . Let $G$ be an Abelian group. If $A \subseteq G$, the subgroup generated by $A$, namely $\< A \> $, is the smallest subgroup of $G$ containing $A$;
if $A$ is the singleton $\{a\}$,
we {just write} $\< a \> $. Tor~$G$ denotes the torsion subgroup of $G$.

{For a \tg \ $G$, we denote by Hom$(G,\TT)$ the set of all homomorphisms $\phi:G\to \TT$. It becomes a group defining the operation pointwise. Then $\wG$, the {\em dual group of $G$,} is defined to be the subgroup of Hom$(G,\TT)$ consisting of continuous functions,  refered to as {\em the characters of} $G$. When $G$ is discrete, its dual group coincides with Hom$(G,\TT)$. 
 A \tg \ $G$ is {\em precompact} if whenever $U$ is an open subset of $G$, there is a {\em finite} subset $F\subseteq G$ \st \ {$G=F+U$}. If in addition to be precompact, $G$ is Hausdorff, then we say that $G$ is {\em \tb .} It is a Theorem of A. Weil \cite{weil1937} that the completion of a precompact (Hausdorff) group $G$ is a compact (Hausdorff) group $\overline{G}$, which we call its {\em Weil completion.}} 
Equipped with the finite-open topology $\gs(\wG,G)$,
$\wG$ becomes a totally bounded \tg ,  which we denote as $\wG_p$. As such, it has a (Weil) completion, referred to as {\em its Bohr compactification} and denoted by $b\wG_p$. {It is important to notice that $\wG_p$ is a dense subgroup of $b\wG_p$, and that the underlying group of $b\wG_p$ is precisely Hom$(G,\TT)$.}
We know $\wZ=\TT$, and $\widehat{\TT}=\ZZ$ by duality, the latter when equipped with the compact-open topology.
If $S$ is a subgroup of Hom$(G,\TT)$, denote by $G_S$ the \tg \ obtained by equipping $G$ with the weakest topology $\tau_S$ that makes the elements of $S$ continuous. It follows that $G_S$ is a precompact group that is Hausdorff \sii \ $S$ separates the elements of $G$ (\ie \ if $S$ is {\em point-separating}) \sii \ $S$ is dense in Hom$(G,\TT)$, which is compact when viewed as $\wG_p$ for $G$ discrete. By the Comfort-Ross Theorem \cite{ComfortRoss1964}, $(\widehat{G_S})_p=S$.
{A topological group $G$ is said to be \emph{Maximally Almost Periodic} ({\it MAP} for short) if $\wG$ separates the points in $G$. {In such a case,} {$\wG$ is dense in $b\wG_p$} 
(see \cite[Theorem 1.12]{ComfortRoss1964}).

{The paper is organized as follows: In Section 2, we develop the algebraic and topological tools necessary to prove the main result, which is proven at the end of the section. In Section 3, we apply the main result to solve a few additional questions and clarify some situations, as mentioned in the introduction.}

\section{Independent sets}
\bdfn\label{def1} 
Let $G$ be an abelian group and let $x\in G$. We define $S(x)=\TT$ if $\text{ord}(x)=\infty$ and $S(x)=\ZZ_q$ if $\text{ord}(x)=q$.

{A subset $A$ of $G$ is said to be \emph{independent} if {$0\not\in A$, and} whenever we have $m_1x_1+\dots +m_nx_n=0,$ { with} $m_j\in\ZZ, x_j\in A, 1\leq j\leq n$, $n\in\NN$,
it holds that $m_jx_j=0, \forall\ 1\leq j\leq n$. By \cite{fuchs} (1.4), $A=\{a_i : i \in I\}$ is independent if and only if $\< A \> = \oplus_{i\in I} \< a_i \> $.}
\edfn

\blem\label{lema4} 
Let $G$  be an abelian group and let $F\in G^{<\omega}$. If $A\subseteq G$ is an infinite independent subset, then there is {an infinite subset} $A'\subseteq A$ \st\ $\langle A' \rangle\cap\langle F \rangle =\{0\}$.
\elem
\pf\ {Since subgroups of finitely generated Abelian  groups are again finitely generated, it follows that $\< F \> \cap \< A\> \leq  \< F \> $ is a finitely generated subgroup of $\< A \> = \oplus_{i\in I} \< a_i \> $. Hence, there exists a finite subset $I_F$ of $ I $ such that $\< F \> \cap \< A\> \subseteq \oplus_{i\in I_F} \< a_i \> $. So $A':= \oplus_{i\in I\setminus I_F} \< a_i \> $ has the required properties.}\epf

\blem\label{lema1} 
Let $G$ be an abelian group and $F\in G^{<\omega}$ an independent subset. If $f\colon F\to \TT$ satisfies that $f(x)\in S(x)\ \forall\ x\in F$,
then there is $\chi\in$ {{\em Hom}$(G,\TT)$} 
\st\ $\chi(x)=f(x)\ \forall\ x\in F$.
\elem
\pf\ Take the subgroup $\langle F\rangle$ generated by $F$ and define $\phi\colon \langle F\rangle\to \TT$ by
$$\phi\left(\sum\limits_{1\leq i\leq n} m_ix_i\right)=\sum\limits_{1\leq i\leq n} m_if(x_i).$$ The map $\phi$ defines a homomorphism since $F$ is an independent subset of $G$.
Furthermore, since $\TT$ is a divisible group, we can extend $\phi$ to a homomorphism $\chi\colon G \to \TT$. {Obviously, $\chi\in$ Hom$(G,\TT)$. \epf}

For the next result, see also \cite{HR} (26.14):

\blem\label{lema2} 
Let $G$ be a MAP abelian group and $F\in G^{<\omega}$ an independent subset. If $f\colon F\to \TT$ satisfies that $f(x)\in S(x)\ \forall\ x\in F$, then for all $\epsilon >0$ there is $\chi\in \widehat{G}$ \st\ $|\chi(x)-f(x)|<\epsilon \ \forall\ x\in F$.
\elem
\pf\ Applying Lemma \ref{lema1}, take $\psi\in b\widehat{G}_p$ \st\ $\psi(x)=f(x)\ \forall\ x\in F$.
Since $\widehat{G}$ is dense in $b\widehat{G}_p$, given $\epsilon >0$, there is $\chi\in \widehat{G}$ \st\
$|\psi(x)-\chi(x)|<\epsilon$ for all $x\in F$.\epf

\blem\label{lema3} 
Let $G$  be a MAP abelian group, $F\in G^{<\omega}$ and $\psi\in b\widehat{G}_p$.
If $g\in G\setminus F$ satisfies that
$\langle g \rangle\cap \langle F \rangle =\{0\}$ and $t\in S(g)$, then for all $\epsilon >0$ there is $\chi\in \widehat{G}$ \st\ $|\chi(x)-\psi(x)|<\epsilon \ \forall\ x\in F$
and $|t-\chi(g)|<\epsilon$.
\elem
\pf\ Consider $\psi_{|\langle F \rangle }\colon \langle F \rangle\to \TT$. Since $\langle g \rangle\cap \langle F \rangle =\{0\}$ and $t\in S(g)$,
we can define $\phi\colon \langle F\cup \{g\} \rangle\to \TT$ by $\phi(x+ng)=\psi(x)+nt$ for all $x\in \langle F \rangle$ {and $n \in \ZZ$}.
We now extend $\phi$ to $G$ as above, using that $\TT$ is divisible. Thus we have that $\phi\in b\widehat{G}_p$.
Again, using the density of $\widehat{G}$ in $b\widehat{G}_p$, we obtain an element $\chi\in\widehat{G}$ such that
$|\chi(x)-\phi(x)|<\epsilon \ \forall\ x\in F\cup \{g\}$. \epf

\blem\label{lema5} 
Let $G$ be a MAP abelian group \st\ $\widehat{G}_p$ is Baire and let $A\subseteq G$ be an infinite independent subset.
Then for any sequence sequence $\{I_k\}_{k<\omega}$ of open subsets in $\TT$ \st\ each $I_k$ contains at least a $n$-root
of the unity for all $2\leq n<\omega$, the set $$N:= \{ \chi\in\widehat{G} : \exists\ \{x_k\}_{k<\omega}\subseteq A : \chi(x_k)\in I_k\ \forall k < \omega\}$$
is a dense $G_\delta$ subset of $\widehat{G}_p$.
\elem
\pf\ {Set $U_k=\bigcup\limits_{a\in A} \{\chi\in\widehat{G} : \chi(a)\in I_k\}$. Clearly $U_k$ is open in $\widehat{G}_p$.
We will show that $U_k$ is also dense in $\widehat{G}_p$. Indeed, let $F\in G^{<\omega}$ and $U$ be an open subset of $\TT$.
Set $N(F,U)=\{\psi\in\widehat{G} : \psi(F)\subseteq U\}$ that we assume nonempty. Pick an arbitrary element $\psi\in N(F,U)${. Since} $F$ is finite, 
there is $\delta >0$ \st\ $(\psi(g)-\delta,\psi(g)+\delta)\subseteq U$ for all $g\in F$ and, since $A$ is an independent subset, by Lemma \ref{lema4},
there is an infinite subset $B \subseteq A$ \st\ $\langle B \rangle\cap\langle F \rangle =\{0\}$. Pick an arbitrary element $b\in B$.
Because $I_k$ contains at least a $n$-root of the unity for all $2\leq n<\omega$,
we can select a point $t\in S(b)\cap I_k$ and a positive real number $0<\epsilon < \delta $ such that $(t-\epsilon, t+\epsilon)\subseteq I_k.$
By {Lemma \ref{lema3}}, 
there is $\chi\in\widehat{G}$ \st\
$|\chi(g)-\psi(g)|<\epsilon<\delta \ \forall\ g\in F$ and $|t-\chi(b)|<\epsilon$, which yields $\chi(b)\in I_k$.
Therefore $\chi\in N(F,U)\cap U_k$, which implies that $U_k$ is dense in $\widehat{G}_p$, which is a Baire space by hypothesis.
Hence $\bigcap\limits_{k<\omega} U_k\not=\emptyset$ is a dense $G_\delta$ subset of $\widehat{G}_p$.
This completes the proof.}\epf
\medskip

To explain the main idea behind the proof of Theorem \ref{teor1} we still need a few results.
{Recall that a space is said to be {\em scattered} if each of its {non-empty} subspaces has an isolated point.  The following is known {(see \cite[8.5.4]{Semadeni} for (1) and \cite[Th. 4]{MRS} for (2)).}

\blem\label{lem1}  
\begin{enumerate}
		\item If $ X $ is an infinite compact space which is not scattered, then there exist surjective continuous mapping
$ f_1 : X \to\TT $. If, in addition, $X$ is $0$-dimensional, then there is a surjective continuous mapping $ f_2 : X \to\{0,1\}^\omega $.
	\item A compact infinite space $K$ that is scattered contains a non-trivial convergent sequence.
\end{enumerate}
\elem
\medskip

The proof of the following result is a variation of the arguments used in \cite[pag. 213]{JNR}. We include it here for completeness' sake.

\blem\label{lem2} 
A compact infinite space $K$ that is not scattered contains a closed subspace $L$ \st \ $L$ has a a countable cover $\{V_n\}$ of open sets, each {containing a compact subspace} of cardinality greater than or equal to $\cc$, \st \ for every open set $V$ of $L$ there is $n$ with $V_n \subseteq V$.\elem

\pf \ {Set $C=\TT$ (respectively, $C=\{0,1\}^\omega$ in case $K$ is $0$-dimensional). }

\noindent {{\em Claim.} There is a closed subspace $L$ of $K$ \st \ there is an {\em irreducible} continuous surjection $f:L\to C$, \ie \ $f$ has the property that if $A\subseteq L$ is closed and proper, then $f[A]\neq C$.} 

{
	{Indeed,} fix a continuous surjective map $f:K \to C$ (which exists by Lemma \ref{lem1} (1)) and let\[\sF:=\{X \subseteq K: X \mbox{ is closed and } f[X]=C\}.\]Plainly, this set is inductively ordered and the maximal element has the desired properties  (See also \cite[Prob. 3.1.C(a)]{engel}).}


We can now prove the Lemma. By the Claim, there {is a} closed subspace $L$ of $K$ \st \ there is an irreducible continuous surjection
$f:L\to {C}$. Let $\{B_n\}$ be a  countable base of ${C}$.
We have that {each $B_n$ contains a compact subspace, say $K_n$, with $|K_n|=\cc$ for all $n \in \NN$. Let $V_n:=f^{-1}[B_n]$ and $C_n:=f^{-1}[K_n]$. Then $C_n$ is compact, $|C_n|\geq\cc$, and $C_n\subseteq V_n$  for all $n \in \NN$.}
Let $V$ be an open subset of $L$.  Then $L\setminus V$ is closed in $L$, hence $F:=f[L\setminus V]$ is closed in ${C}$.
Since $f$ is irreducible, $B:=C \setminus F$ is non-empty, and obviously open, hence there is a $B_n$ contained in $B$.
It follows that $V_n \subseteq f^{-1}[B] \subseteq V$, as required.}
\epf
\medskip

The following lemma is crucial in order to demonstrate our main result. Our proof uses an argument that appears in \cite[pag. 213]{JNR} in which it is attributed to ``folklore". We are very grateful to Professor Roman Pol for pointing out this article to us
{(cf. \cite{MykPol}).}

\blem\label{lema6.1} 
Let $G$ be a totally bounded abelian group \st\ every uncountable compact subset of $G$
contains an uncountable independent subset. If the dual group $\widehat{G}_p$, equipped with
the finite-open topology, is a Baire group, then every compact subset of $G$ is finite.
\elem
\pf\
Let $K$ be an infinite compact subspace of $G$. If $K$ {were scattered, by {Lemma \ref{lem1} (2), it would contain a convergent sequence, which is impossible by} Corollary 2.4 in \cite{BruTka2012}. We can therefore assume that $K$ is not scattered.

We take $\TT:=(-1/2,1/2], I_{2m}:=(-1/8,1/8)$, and $I_{2m+1}=(-1/2,-1/4) \cup (1/4,1/2]$ for all $m<\omega$. Consider the closed subset $L$ of $K$, {with} the countable collection of open subsets $\{V_n\}$, {each containing a compact subset of cardinality $ \geq \cc$, provided by Lemma \ref{lem2}. Thus, according to the hypotheses we have assumed, each $V_n$ contains an uncountable independent subset, say $A_n$. Define
\[S_n:=\{f \in \widehat{G}_p: \exists\ \{x_{nk}\}_{k<\omega}\subseteq A_n\ \hbox{with}\ f(x_{nk})\in I_{k}\ \forall k<\omega\}.\]}
By Lemma \ref{lema5}, we have that $S_n$ contains a dense $G_\delta$ subset of $\widehat{G}_p$ for all $n<\omega$.

Since $\widehat{G}_p$ is Baire by hypothesis, we have that $\bigcap_{n=1}^\infty S_n$ contains a dense $G_\delta$ subset of $\widehat{G}_p$.
Now, pick $\phi \in \bigcap_{n=1}^\infty S_n$. Let $w\in L$ and consider $U=\{t \in \TT: |t-\phi(w)|<1/16\}$. By continuity, $V:=\phi^{-1}[U]$ is a \nhd \ of $w$ in $G$
and, by Lemma \ref{lem2}, there is a $V_n$ inside $V \cap L$.
Since $\phi\in S_n$, there exists $\{x_{nk}\}_{k<\omega}\subseteq V_n$ with $\phi(x_{nk})\in I_k\ \forall k<\omega$. But then
\[1/8< |\phi(x_{n0})-\phi(x_{n1})|\leq |\phi(x_{n0})-\phi(w)|+|\phi(w)-\phi(x_{n1})|<1/8,\]
a contradiction. This completes the proof. \epf
\medskip

Let $G$ be an abelian group and let $m$ be a natural number. In what follows we write $G(m):=\{g\in G : mg=0\}$.
If $H$ is a subset of $G$, we write $H^\bot:=\{\phi\in \wG: \phi[H]=\{0\}\}.$
Lemmata \ref{lema5.1} and \ref{lema5.2} follow from the main result in \cite{ChascoDominguezTkachenko}, which make use
of quite involved techniques. However, we include the proofs here in order to make the paper self-contained. In the sequel, we say that an abelian group $G$ has order $m\in\NN$ if $m$ is the smallest
number satisfying $mg=0$ for all $g\in G$.

\blem\label{lema5.1} 
Let $G$ be a bounded torsionMAP abelian group, whose order is a prime number $p$, \st\ $\widehat{G}_p$ is Baire.
Then every compact subset of $G$ must be finite.
\elem
\pf\
First, observe that $G$ is a vector space on the field $\ZZ(p)$. Let $K$ be an uncountable subset of $G$ and
let $F$ be a maximal independent subset of $K$. If we assume that $F$ is countable, the subgroup $\langle F\rangle$, generated by $F$,
is also countable. Since $G$ is a vector space on $\ZZ(p)$, it follows that for every element $x\in K\setminus \langle F\rangle$,
the subset $F\cup\{x\}$ is linearly independent. This fact contradicts the maximality of $F$. 
Thus, every uncountable compact subset $K$ of $G$ must contain an uncountable maximal independent subset.
Hence, we are in a position to apply Lemma \ref{lema6.1}, which completes the proof.
\epf

\blem\label{lema5.2} 
Let $G$ be a bounded torsionMAP abelian group \st\ $\widehat{G}_p$ is Baire.
Then every compact subset of $G$ must be finite.
\elem
\pf\
We apply induction on the order of $G$. We have already verified the result when the order of $G$ is a prime number.
Now, assume that the result holds for every group of order less than $m\in\NN$, and let $G$ be a group of order $m$.
Again, {observe} that $G$ does not contain {infinite countable} compact sub\-sets since the group $\widehat{G}_p$ is Baire
(see \cite[Corollary 2.4]{BruTka2012}). Thus, reasoning by contradiction, suppose that that $K$ is an uncountable
compact subset of $G$. By Lemma \ref{lema6.1}, in order to obtain a contradiction,
it will suffice to verify that $K$ contains an uncountable independent subset.
Again, reasoning by contradiction, let $F$ be an independent maximal subset of $K$ and
sup\-pose that $F=\{x_n\}_{n<\omega}$ is finite or countably infinite.
Then for every $x\in K$ there are $m_x\in\ZZ$ and $\overline{m}_x\in\bigoplus\limits_{n<\omega}\ZZ$ \st\
$$m_xx=\overline{m}_x(n_1)x_{n_1}+\dots + \overline{m}_x(n_k)x_{n_k}.$$
There is no loss of generality in assuming that $m_x$ is the smallest natural number satisfying the above equality.
Since both sets $F^{<\omega}$ and $\bigoplus\limits_{n<\omega}\ZZ$ are countable, there must be
$\overline{m}\in\bigoplus\limits_{n<\omega}\ZZ$, $\{x_{n_1},\dots , x_{n_k}\}\in F^{<\omega}$ and an uncountable infinite
 compact subset $\widetilde{K}\subseteq K$ \st\
$$m_xx=\overline{m}(n_1)x_{n_1}+\dots + \overline{m}(n_k)x_{n_k}\ \forall x\in \widetilde{K}.$$
Furthermore, since $\ZZ$ is countable, there must be $m_0\in\ZZ$ and
an uncountable infinite compact subset $\widetilde{\widetilde{K}}\subseteq \widetilde{K}$ \st\
$$m_0x=\overline{m}(n_1)x_{n_1}+\dots + \overline{m}(n_k)x_{n_k}\ \forall x\in \widetilde{\widetilde{K}}.$$
Observe that if $m_0=m$, then
$$0=m_0x=\overline{m}(n_1)x_{n_1}+\dots + \overline{m}(n_k)x_{n_k}\ \forall x\in \widetilde{\widetilde{K}}.$$
This means that $$\langle x\rangle\cap \langle F\rangle =\{0\}.$$
Thus $F\cup\{x\}$ is an independent set, in contradiction with the maximality of $F$.
Therefore, we may assume that $m_0<m$ without loss of generality.

If we fix a point $x_0\in \widetilde{\widetilde{K}}$ and set $L=\widetilde{\widetilde{K}}-{x_0}$,
we obtain an uncountable infinite compact subset of $G$ \st\ $m_0x=0$ for all $x\in L$.
That is $L\subseteq G(m_0)=\{g\in G : m_0g=0\}$. The group $\widehat{G}/G(m_0)^\bot$
equipped with the topology $t_p(G(m_0))$ is a continuous homomorphic image of the group $\widehat{G}_p$,
which is a Baire group by hypothesis. Therefore, the group $(\widehat{G}/G(m_0)^\bot,t_p(G(m_0)))$ is a Baire group,
and its dual group $G(m_0)$ contains $L$, has order $m_0$, and contains an uncountable compact subset. This is impossible
according to our inductive hypothesis. Therefore we have reached a contradiction that completes the proof.
\epf
\medskip

We can now formulate our main result that 
solves positively  Question 3.4 in \cite{ChascoDominguezTkachenko}.
This result sharpens Lemma \ref{lema5.2} with a similar proof.

\bthm\label{teor1} 
Let $G$ be a totally bounded abelian group \st\ its dual group $\widehat{G}_p$, equipped with
the finite-open topology is a Baire group. Then every compact subset of $G$ must be finite.

\ethm
\pf\  {Again, observe that $G$ does not contain {infinite countable} compact sub\-sets since the group $\widehat{G}_p$ is Baire
(see \cite[Corollary 2.4]{BruTka2012}). {Let $K$ be an uncountable compact subspace of $G$. By Lemma \ref{lema6.1} we must show that $K$ contains an uncountable independent set.}}
Thus, reaso\-ning by contradiction, let $F$ be an independent maximal subset of $K$ and
sup\-pose that $F=\{x_n\}_{n<\omega}$ is finite or countably infinite. 
Then for every $x\in K$ there are $m_x\in\ZZ$ and $\overline{m}_x\in\bigoplus\limits_{n<\omega}\ZZ$ \st\
$$m_xx=\overline{m}_x(n_1)x_{n_1}+\dots + \overline{m}_x(n_k)x_{n_k}.$$
Since both sets $F^{<\omega}$ and $\bigoplus\limits_{n<\omega}\ZZ$ are countable, there must be
$\overline{m}\in\bigoplus\limits_{n<\omega}\ZZ$, $\{x_{n_1},\dots , x_{n_k}\}\in F^{<\omega}$ and an uncountable infinite
 compact subset $\widetilde{K}\subseteq K$ \st\
$$m_xx=\overline{m}(n_1)x_{n_1}+\dots + \overline{m}(n_k)x_{n_k}\ \forall x\in \widetilde{K}.$$
Furthermore, since $\ZZ$ is countable, there must be $m_0\in\ZZ$ and
an uncountable infinite compact subset $\widetilde{\widetilde{K}}\subseteq \widetilde{K}$ \st\
$$m_0x=\overline{m}(n_1)x_{n_1}+\dots + \overline{m}(n_k)x_{n_k}\ \forall x\in \widetilde{\widetilde{K}}.$$
If we fix a point $x_0\in \widetilde{\widetilde{K}}$ and set $L=\widetilde{\widetilde{K}}-{x_0}$,
we obtain an uncountable infinite compact subset of $G$ \st\ $m_0x=0$ for all $x\in L$.
That is $L\subseteq G(m_0)=\{g\in G : m_0g=0\}$. The group $\widehat{G}/G(m_0)^\bot$
equipped with the topology $t_p(G(m_0))$ is a continuous homomorphic image of the group $\widehat{G}_p$,
which is a Baire group by hypothesis. Therefore, the group $(\widehat{G}/G(m_0)^\bot,t_p(G(m_0)))$ is a Baire group
and its dual group $G(m_0)$ has order $m_0<m$ and contains $L$, an uncountable compact subset. This is impossible
according to Lemma \ref{lema5.2}. Therefore we have reached a contradiction that completes the proof. \epf

\section{Some applications of the main result}

In \cite{ChMPVa:99} the following notion is defined: a topological abelian group $G$ is said to be \emph{$g$-barrelled} if any compact subset of $\widehat{G}_p$ is equicontinuous. Among other results, the authors
prove that hereditarily Baire groups or \v{C}ech-complete groups are $g$-barrelled. As a consequence of Theorem \ref{teor1}, we identify a new class of $g$-barrelled groups.

\bcor\label{cor21} 
Every totally bounded Baire group is $g$-barrelled.
\ecor
{\pf \ Assume that $G$ is a totally bounded Baire group. By the Comfort-Ross Theorem \cite{ComfortRoss1964} or Theorem 1 of \cite{RaczTrig}, we have that $G$ and $\widehat{(\wG_p)}_p$ are topologically isomorphic. Therefore, the dual $\widehat{(\wG_p)}_p$ of $\widehat{G}_p$ is a totally bounded Baire group.
By Theorem \ref{teor1}, every compact subspace of $\widehat{G}_p$ is finite, hence equicontinuous, as required. \epf}

\brem\label{cor21lqcvx} One {cannot} relax the ``\tb" hypothesis. 
See Theorem 6.4.1 of \cite{AusDik:2021}.
\erem

The notion of groups in duality was introduced by Varopoulos in \cite{varopoulos:64}
for topological Abelian  groups as follows:

\bdfn[Varopoulos]\label{Varo1} 
Let $G$ and $G'$ be two Abelian topological groups. Then we say
that they are in duality if and only if there is a function
$$\langle\cdot,\cdot\rangle:G\times G'\rightarrow \TT\quad \text{ such that }$$

\begin{itemize}
\item [a)]$\langle g_1 + g_2,g'\rangle=\langle g_1,g'\rangle+
\langle g_2,g'\rangle$ for all $g_1,g_2\in G$, $g'\in G'$;
\item [b)]$\langle g,g'_1+ g'_2\rangle=\langle g,g'_1\rangle+
\langle g,g'_2\rangle$ for all $g\in G$, $g'_1,g'_2\in G'$;
\end{itemize}
and the following hold:
\begin{itemize}
\item [i)] if $g\not= 0_G$, the neutral element of $G$,  then there
exists $g'\in G'$ such that
$\langle g,g'\rangle\not= 0$;
\item [ii)] if $g'\not= 0_{G'}$, the neutral element of $G'$, then
there exists $g\in G$ such that
$\langle g,g'\rangle\not= 0$.
\end{itemize}
\edfn

\bdfn\label{Varo2} 
If we have a duality $\langle G,G'\rangle$, we say that a topology $\tau$ on
$G$ is \emph{compatible} with the duality when {when the groups $(G,\tau)\ \widehat{ }$  and $G'$ are algebraically isomorphic.}
A locally quasi-convex Hausdorff group $(G,\tau)$ is called a \emph{Mackey group} if {$\tau$ is
the supremum of all the compatible topologies with the duality $(G,\widehat{G})$,
and $\tau$ itself  is compatible with this duality (see \cite{ChMPVa:99}).}
\edfn

It is known (ibidem) that every $g$-barrelled locally quasi-convex Hausdorff group is a Mackey group.
Since from Corollary \ref{cor21} we know that every totally bounded Baire group is $g$-barrelled, we obtain:

\bcor\label{cor23}  
Every totally bounded Baire group is a Mackey group.
\ecor
\mkp

{If $X$ is a Tychonoff space, 
and $A(X)$ denotes the (Markoff) free Abelian \tg \ on $X$, then its character group can be identified with the group of continuous functions from $X$ to $\TT$, which we denote as $C(X,\TT)$. 
{Recall that if $G$ is a \tg, the {\em Bohr topology} of $G$ is the weakest topology on $G$ that makes the elements of $\wG$ continuous. When $\wG$ separates its points, $G$ equipped with its Bohr topology is a totally bounded group. We denote by $A_p(X)$ the
group $A(X)$ equipped with its Bohr topology.
Writing $G=A_p(X)$, we obtain $\widehat{G}_p=C_p(X,\TT)$, \ie \ $C(X,\TT)$ equipped with the finite-open topology
(see \cite{gh99}). By Theorem \ref{teor1}, we have that if $G$ has an infinite}
compact subspace, then $\widehat{G}_p$ can't be Baire. Applying this fact to the group $A(X)$, we obtain:}

\bcor \label{cor1} 
{If $X$ is a Tychonoff space containing an infinite compact subspace, then $C_p(X,\TT)$ is not Baire. In particular, this is the case if $X$ is compact and infinite.}
\ecor

{\pf \ If $K$ is an infinite compact subspace of $X$, then $K$ is an infinite compact subspace of $A(X)$, hence of $G=A_p(X)$, therefore $\wG_p = C_p(X,\TT)$ cannot be Baire. \epf}


\bcor \label{cor2} 
If $K$ is a compact zero-dimensional space that is not scattered, then $C_p(K,\{0,1\})$, is not Baire.
\ecor

{\pf \ {By Lemma \ref{lem2}, $K$ contains a compact set $ L $ which has a continuous irreducible surjection $ f : L \to \{0,1\}^\omega $. Consider the sets $V_n$ provided by Lemma \ref{lem2}. The set $ G_n $ of functions in $ C_p(L, \{0,1\}) $ which take both values 0 and 1 on $ V_n $ is open and dense in $C_p(L, \{0,1\})$: For, because $L$ is Tychonoff and each $V_n$ is infinite, it follows that $G_n$ is not empty. Suppose $f\in G_n$. Then there is $x \in V_n$ and $y\in V_n$ \st \ $f(x)=0$ and $f(y)=1$. Then $X:=\{g\in C_p(L, \{0,1\}):g(x)=0\}$ and $Y:=\{g\in C_p(L, \{0,1\}):g(y)=1\}$ are both open, and $f\in X\cap Y \subseteq G_n$. Hence, each $G_n$ is open. To see that it is dense, let $F_0$ and $F_1$ be finite subsets of $L$ and consider the open set $V:=\{f\in C_p(L, \{0,1\}): f[F_j]=\{j\}, j=0,1\}$. Pick $x_0 \in V_n$ and a distinct $x_1\in V_n$. Using that $L$ is Tychonoff find a continuous function $g:L\to \{0,1\}$ with $f[F_j\cup \{x_j\}]=\{j\}$. Then, $f\in G_n \cap V$, as required. Any function in the intersection of all $ G_n $ would have  oscillation 1 at each point of $ L $ so the intersection is empty.\epf}
\mkp

Professor Roman Pol communicated this result and its proof to us in a private email.

\brem \label{Mykh} In particular, Theorem 4.1 in \cite{Myk} considers the case $Y=\beta \omega \setminus\omega$ and shows that the following are equivalent for $C_p(Y,\{0,1\})$:
\begin{enumerate}
	 \item being meagre,
	 \item not being Baire,
	 \item the existence of a sequence of finite pairwise disjoint sets $E_n \subseteq Y$ converging {\em weakly} to one point in $Y$,
	 \item the existence of a sequence of finite pairwise disjoint sets $E_n \subseteq Y$ converging {\em weakly} to every point in $\cup E_n$.
\end{enumerate}
By Corollary \ref{cor2}, $C_p(\beta \omega \setminus\omega,\{0,1\})$ satisfies all the above conditions.

\erem

{Let $X$ be a Tychonoff space.} Since $e:[0,1]\to \TT$, with $e(x)=x$ for $x\in[0,1)$ and $e(1)=0$ is continuous, it follows from \cite[Problem 91]{Tk2011} that the map $h_e:C_p(X,[0,1]) \to C_p(X,\TT)$ defined by $h_e(f)=e\circ f$ is continuous.} 



\bthm \label{teor2} 
{Let $X$ be a Tychonoff space.} If $C_p(X,[0,1])$ is Baire, then $C_p(X,\TT)$ is also Baire.
\ethm

\pf {\ Let $ C_{01} = \{f\in C_p(X, [0, 1]): \{0,1\} \cap f[X] = \emptyset\}$ and $C_1 = \{g \in C_p(X,\TT): 1 \not\in g[X]\} $. For every open subset
	$ V $ of $ C_p(X, [0, 1]) $, the set $ V \cap C_{01} $ is dense in $ V $. And for every open set $ U $ of $C_p(X,\TT)$, the set $ U \cap C_1 $ is dense in $ U $. Further,
	$h_{e|C_{01} }$ induces a homeomorphism of $ C_{01} $ to $ C_1 $.}

{If $ (U_n) $ are open, dense subsets in $ C_p(X,\TT) $, then it follows easily that
$ (h^{-1}_e [U_n]) $ are open dense subsets in $ C_p(X, [0,1]) $. Since $ C_p(X, [0,1]) $ is assumed to be Baire, the set $\bigcap h^{-1}_e [U_n] = h^{-1}_e \bigcap [U_n]$ is dense in $ C_p(X, [0,1]) $. Since $ h_e $ is continuous, $h_e(h^{-1}_e [\bigcap [U_n]])$ is dense in $ h_e[ C_p(X, [0,1]] $ and hence $\bigcap [U_n]$ is dense in $ C_p(X,\TT) $.\epf}
}
\mkp

{We thank the refereee for suggesting this short proof of Theorem \ref{teor2}.}

\bcor\label{BaireNoinfinitecpt} 
{Let $X$ be a Tychonoff space.} If $C_p(X,[0,1])$ is Baire, then $X$ cannot have infinite compact subsets.
\ecor
{\pf \ By Theorem \ref{teor2}, it follows that $\widehat{A_p(X)}_p=C_p(X,\TT)$ is also Baire. Therefore $A_p(X)$ cannot have infinite compact sets by Theorem \ref{teor1}. If $K$ were compact with $K \subseteq X\subseteq A(X)$, then $K \subseteq A_p(X)$ would also be compact, hence it should be finite. \epf}

\brem \label{Tkachuk} 
In Problem 286 in \cite{Tk2011}, it is constructed an infinite dense pseudocompact $X \subseteq [0,1]^\cc$ \st \ each of its countable subspaces is closed and discrete, hence this space cannot have infinite compact subspaces. It also shown (loc. cit.) that $C_p(X,[0,1])$ is Baire, despite the fact that $C_p(X)$  is not. This and Corollary \ref{BaireNoinfinitecpt} provide another proof that $X$ cannot have infinite compact subspaces.
\erem


The following definition can be found in \cite[Def. 5.4.1]{AusDik:2021}.

\bdfn\label{AusDik:2021} 
Let $K$ be a topological group and $u = (u_n$) be a sequence in $\widehat{K}$.
Let $s_u(K) =\{ x\in K : \lim_n u_n(x) = 0\ \hbox{in}\ \TT\}$. A subgroup $G$ of $K$ is called
\begin{itemize}
  \item \emph{characterizable} (or \emph{characterized}), if there exists $u = (u_n)$ in $\widehat{K}$ such that $G = s_u(K)$;
  \item \emph{$g$-closed}, if $G$ is the intersection of characterizable subgroups of $K$;
  \item \emph{$g$-dense}, if $G\subseteq s_u(K)$ for some $u = (u_n)$ in $\widehat{K}$ yields $s_u(K) = K$.
\end{itemize}
\edfn
\mkp

In connection with this definition, Au\ss enhofer and Dikranjan ask whether every totally bounded group $G$ that is
$g$-dense in its completion must be a $g$-barrelled group (see \cite[Question 11.3.3]{AusDik:2021}). Next we give a negative answer to this question.

\bcor\label{cor24} 
Set $\omega^*:=\beta\omega\setminus \omega$. Then the group $C_p(\omega^*,\TT)$ is $g$-dense in its completion but not $g$-barrelled.
\ecor
\pf\ According to \cite[Theor. 4.19]{gh:FM}, we have that $A(\omega^*)$ strongly respects
compactness. This means that every compact subset $C$ of $A(\omega^*)_p$ is also compact in the free topological abelian group $A(\omega^*)$,
which implies that there is $n<\omega$ \st\ $C\subseteq \langle\omega^*\rangle_n$ {\ie \ words of length $\leq n$}. As a consequence, $C$ {cannot} contain nontrivial convergent
sequences. According to Lemma \cite[5.4.5]{AusDik:2021}, this means that the group  $C_p(\omega^*,\TT)$ is $g$-dense in its completion.
On the other hand, its dual group contains $\omega^*$, an infinite compact subset. Thus, {being totally bounded,
the group $C_p(\omega^*,\TT)$ cannot be} $g$-barrelled \cite[5.1.24]{AusDik:2021}.\epf
\medskip


Examples of totally bounded Pontryagin-reflexive groups have been provided by different authors. Here we collect some of them:
It is shown in {\cite[Theorem 2.8]{ArdChaDomTka:2012} and \cite[Theorem 6.1]{gm2011},} that a pseudocompact group without infinite compact subsets is
reflexive. A slightly more general fact is established in {\cite[Theorem 2.8]{BruTka2012}}: Every totally bounded Baire group without infinite compact subsets is reflexive
provided that it satisfies the so-called Open Refinement Condition. More recently, it has been proved in \cite{ChascoDominguezTkachenko}
that if  $ G $ is a bounded torsion totally bounded group which is a Baire space without infinite compact subsets, then $ G $ is Pontryagin reflexive.
This result has been subsequentlty extended to torsion totally bounded groups which are Baire spaces without infinite compact subsets \cite{AusDik:2021}.
The following Corollary extends all these results.

\bcor\label{cor3.15} 
Every totally bounded Baire group without infinite compact subsets is reflexive.
\ecor
\pf\ {Let $G$ be a totally bounded Baire group. 
Recall that a \tg \ $G$ is reflexive if the evaluation map $\Omega:G \to \widehat{(\wG_k)}_k $ is a topological isomorphism, where the subindex $k$ denotes the compact open topology. Since $G$ contains no infinite compact subsets, we have that $\wG_k=\wG_p$. 
 Since $G$ is Baire, it follows that $\wG_k=\wG_p$ has no infinite compact sets, by Theorem \ref{teor1}. Hence, $\widehat{(\wG_k)}_k=\widehat{(\wG_p)}_k=\widehat{(\wG_p)}_p$. By the Comfort-Ross Theorem \cite{ComfortRoss1964} or Theorem 1 of \cite{RaczTrig}, the evaluation map $\Omega:G \to \widehat{(\wG_p)}_p $ is a topological isomorphism. But $\widehat{(\wG_p)}_p=\widehat{(\wG_k)}_k$, as required}.\epf

\brem \label{TkaProb3.9} 
In \cite{ChascoDominguezTkachenko} it is also asked (Problem 3.9) the following. Let $G$ be a precompact Boolean group \st \ its dual $\wG_p$ has no infinite compact sets. Does $G$ have the Baire property?
As shown in Theorems 2.3 and 3.9 of \cite{FerHerTka} and Example 5.3.13 of \cite{AusDik:2021}, the answer is ``not always". Here we offer yet another counterexample. By Theorem 1.9 in \cite{hartku06}, if $X$ is an infinite compact Abelian {\em metrizable} group, then it contains a dense subgroup $S$ of measure 0 \st \ $(\widehat{X},t_p{(S)})$  has no infinite compact sets.
In the same paper it is proven (Proposition 5.3) that such $S$ is a Borel set. Consider $G:= \bigoplus_\omega \{0,1\}$. Then $X:=\wG=\{0,1\}^\omega $ is a compact metric group  containing a dense Borel subgroup $S$ of measure 0 \st \  $(\widehat{X},t_p{(S)})$  has no infinite compact sets. We claim that $S$ is not of second category. If it were, then $S$ would be nonmeager by \cite{sriva98} (3.5.7). {Furthermore, by} \cite{sriva98} (3.5.13), $S$ would be clopen in $\wG$, and by density, $S=\wG$, contradicting that $mS=0$.\epf \erem
\medskip

{It is easy to see that if $G$ is a \lc \ Abelian non-compact group, then $(G,t_p(\wG))$ is meager. 
We propose the following:}

\bprb \label{G+meager} 
Characterize those MAP groups $G$ \st \ $(G,t_p(\wG))$ is meager.\eprb

{\section{Acknowledgement}}

{ We are are very grateful to Professor Roman Pol for sharing with us several helpful results.
We also thank the referee for an informative and constructive report and to Lydia Au\ss enhofer and Elena Mart\'\i n-Peinador for incisive readings of an early version of this article.}


\begin{thebibliography}{10}
\bibitem{ArdChaDomTka:2012} {Ardanza-Trevijano, S., Chasco, M. J., Domínguez, X. and Tkachenko, M. G., {\em Precompact non-compact reflexive Abelian groups}, Forum Math. 24 (2) (2012)
	289–302.}
	
\bibitem{AusDik:2021}
Au\ss enhofer, L. and Dikranjan, D., \emph{Mackey groups and Mackey topologies}, Dissertationes Mathematicae 567 (2021), 1--141.

	\bibitem{BruTka2012}
	Bruguera, M. and Tkachenko, M., {\em Pontryagin duality in the class of precompact Abelian groups and the Baire property,} Journal of Pure and Applied Algebra, 216 (2012) 2636–2647.
	
	\bibitem{ChascoDominguezTkachenko}
	Chasco, M. J., Dom\'\i nguez, X. and Tkachenko, M., {\em Duality properties of bounded torsion topological Abelian groups,} J. Math. Anal. Appl. 448 (2017), 968-981.

\bibitem{ChMPVa:99}
Chasco, M. J., Mart\'in Peinador, E. and Tarieladze, V., \emph{On Mackey Topology for groups},
Stud. Math. 132, No.3, (1999) 257-284.	
	
\bibitem{ComfortRoss1964}
Comfort, W. W. and Ross, K. A., {\em Topologies induced by groups of characters,} {Fundamenta Mathematicae}, {55}, (1964), {283-291}

	\bibitem{engel}
	Engelking. R., {\em General Topology,} Heldermann Verlag, Berlin 1989.

\bibitem{FerHerTka} Ferrer, M. V., Hern\'andez, S. and Tkachenko, M., {\em On convergent sequences in dual groups}, RACSAM 114, 71 (2020).

{\bibitem{fuchs}
Fuchs, L., {\em Infinite Abelian Groups, vol. 1.} Academic Press, New York, 1970.}

\bibitem{gh:FM} Galindo, J. and Hern\'andez, S., {\em The concept of boundedness and the Bohr compactification of a MAP Abelian group}, Fund. Math. 159 (1999), no. 3, 195218. 159 (1999)
	
	\bibitem{gh99}
	Galindo, J. and Hern\'andez, S., {\em Pontryagin-van Kampen reflexivity for free Abelian \tg s,} Forum Math. 11 (1999), 399-415.
	
\bibitem{gm2011} {Galindo, J. and Macario, S., {\em Pseudocompact group topologies with no infinite compact subsets}, J. Pure Appl. Algebra 215 (4) (2011) 655–663.}

\bibitem{hartku06}
Hart, J. E. and Kunen, K., {\em Limits in compact Abelian groups}, Topology and its Applications 153 (2006), 991–1002.

\bibitem{JNR}
Jayne, J. E., Namioka, I. and Rogers, C. A., {\em Fragmentability and $\sigma$-fragmentability}, Fundamenta Mathematicae,
143 (1993), 207--220.
	
	\bibitem{HR}
	Hewitt, E. and Ross K. A., {\em Abstract Harmonic Analysis: Volume I, Structure of Topological Groups,
		Integration Theory, Group Representations,} Grundlehren der mathematischen Wissenschaften,
		Springer New York, 1963.

\bibitem{MRS}
Mrowka, S., Rajagopalan, M. and Soundararajan, T., {\em A characterization
of compact scattered spaces through chain limits}, TOPO 72 -- General
Topology and its Applications, Lecture Notes in Mathematics 378 (1974), 288--297.

\bibitem{Myk}
Mykhaylyuk, V. V., {\em On questions which are connected with Talagrand problem}, Matematyczni
Studii 29 (2008), 81-88. (in Ukrainian; English version: arXiv:1601.03163 [math.GN], 13
	January 2016).

\bibitem{MykPol}
Mykhaylyuk V. and Pol R., \emph{On a problem of {T}alagrand concerning separately continuous functions},
J. Inst. Math. Jussieu, 20, no. 5, (2021), 1719--1728.

{\bibitem{RaczTrig}
Raczkowski, S. U. and Trigos-Arrieta, F. J., \emph{Duality of totally bounded Abelian groups}, Bolet\'\i n de la Sociedad Matem\'atica Mexicana, 3 (7) (2001), pp. 1-12. MR: 2002g:22006. Zbl. 1007.22004.}

\bibitem{Semadeni}
Semadeni, Z., {\em Banach Spaces of Continuous Functions}, Vol. 1, PWN, Warszawa, 1971.

\bibitem{sriva98}
Srivastava, S. M., {\em A Course on Borel Sets}, Graduate Texts in Mathematics, vol. 180, vol. 180, New York-Berlin-Heildelberg, 1998.

\bibitem{Tk2011}
Tkachuk, V. V., {\em A Cp-Theory Problem Book, Topological and Function Spaces,} Springer, New York, Dordrecht, Heidelberg, London, 2011.

\bibitem{varopoulos:64}
Varopoulos, N., \emph{Studies in harmonic analysis}, Proc. Cambridge Philos. Soc. 60 (1964), 465–516.

\bibitem{weil1937}
Weil, A., \emph{Sur les {E}spaces \`a {S}tr\`ucture {U}niforme et sur la {T}opologie {G}\'en\'erale},
Publ. Math. Univ. Strasbourg, Hermann, Paris, 1937.
	
\end{thebibliography}
\end{document}